\documentclass[a4paper, 11pt, twoside, leqno]{amsart}

\usepackage{styles/generalstyle}
\usepackage{styles/mytheorems}
\usepackage{styles/mystyle}

\title{On the real Section Conjecture in \'etale homotopy theory}
\author{Tim Holzschuh}
\date{\today}

\bibliography{literature.bib}

\renewenvironment{abstract}
 {\par\noindent    \textbf{\abstractname .}  \ \ignorespaces}
 {\par\medskip}

\begin{document}

  \maketitle

  \begin{abstract}
    We study the Section Conjecture in \'etale homotopy theory for varieties over $\rR$.
    We prove its pro-$2$ variant for \emph{equivariantly triangulable} varieties.
    Examples include all smooth varieties as well as all (possibly singular) affine/projective varieties.
    Building on this, we derive the real Section Conjecture in the geometrically \'etale nilpotent (e.g. simply connected) case.	
  \end{abstract}

  \tableofcontents

  \setcounter{section}{0}

  \vspace{-1cm}

  \section{Introduction}
  \label{sec:introduction}
  We study the \emph{Section Conjecture in \'etale homotopy theory} over the real numbers $\rR$.
Let us fix the following situation:

\begin{nul}
    \label{setup:intro}
    Let $k$ be a field with separable closure $\algclos{k} \supset k$, and let $X$ be a geometrically connected quasi-compact and quasi-separated scheme over $k$.
    Write $\Gal_{k} \defeqq \Gal(\algclos{k} / k)$ for the absolute Galois group of $k$, denote by $X_{\algclos{k}}$ the base change of $X$ to $\algclos{k}$, and let $\bar{x} \in X_{\algclos{k}}$ be a geometric point.
\end{nul}

Recall that by \cite[Théorème IX.6.1]{SGA1}, we have the \emph{fundamental exact sequence} 
\begin{equation*}\label{eq:fundamental-exact-sequence}
    \begin{tikzcd}[cramped, sep=small]
        1 \arrow[r]
          & \ethtpygrp[1](X_{\algclos{k}}, \bar{x}) \arrow[r]
            & \ethtpygrp[1](X, \bar{x}) \arrow[r]
              & \Gal_{k} \arrow[r]
                & 1
    \end{tikzcd}   
\end{equation*}
attached to $X/k$.

In his seminal letter \cite{GanF} from 1983 to Faltings, Grothendieck suggested to study the geometry and arithmetic of $X/k$ in terms of its fundamental exact sequence.

Among other things, he conjectured the following.

\begin{sectconj}
  \label{conj:section-conjecture}
    Let $X$ be a connected and proper \emph{hyperbolic curve} over some finitely generated field extension $k / \! \qQ$.
    Then the canonical \emph{Kummer} map
    \[
        \kappa_{X/k} \from X(k) \to \Hom_{\Gal_k}^{\out}(\Gal_k, \ethtpygrp_{1}(X))
    \]
    induced by the functoriality of $\ethtpygrp_{1}(\blank)$ is a bijection.
\end{sectconj}

Except for $k = \cC$, where it is obviously wrong, the analogous statement is also conjectured to hold over local fields $k$ of characteristic $0$.
In this article, we study a generalisation of the above conjecture in terms of \'etale homotopy theory in the case $k = \rR$ of the real numbers.

\subsection{The real Section Conjecture}

Over the reals, injectivity of the Kummer map ceases to hold. 
Instead, two real points $a$ and $b$ of $X$ determine the same conjugacy-class of sections if and only if they lie in the same connected component of the \emph{real analytification} of $X$, henceforth simply denoted by $X(\rR)$.
Mochizuki showed in \cite[Thms. 3.13-3.15]{mzki03} that $\kappa_{X/\! \rR}$ is surjective, even for curves of genus $1$.
Wickelgren later showed in \cite[Thm. 1.1]{Wickelgren}, that one may furthermore replace $\etfdtlgrp(X)$ by its \emph{geometrically pro-$2$ quotient}
\[
  \pquotient[\ethtpygrp_1(X)][(2)] \defeqq \etfdtlgrp(X) \amalg_{\etfdtlgrp(X_{\cC})} \pquotient[\etfdtlgrp(X_{\cC})][2],
\]
that is obtained by pushing out the fundamental exact sequence along the maximal pro-$2$ quotient $\etfdtlgrp(X_{\cC}) \twoheadrightarrow \pquotient[\etfdtlgrp(X_{\cC})][2]$.
We can summarise the state-of-the-art of the real Section Conjecture as follows.

\begin{theorem}[(pro-$2$) real Section Conjecture]
    \label{thm:mochizuki-real-section-conjecture}
    Let $X/\! \rR$ be a smooth, geometrically connected \emph{curve} of arithmetic genus $g \geq 1$.
    All of the maps
    \[
        \htpygrp_{0}\! X(\rR) \to \Hom_{\Gal_{\rR}}^{\out}(\Gal_{\rR}, \ethtpygrp_{1}(X)) \to \Hom_{\Gal_{\rR}}^{\out}(\Gal_{\rR}, \pquotient[\ethtpygrp_1(X)][(2)])
    \]
    are bijective.
\end{theorem}

Our main results are generalisations of \Cref{thm:mochizuki-real-section-conjecture} to higher-dimensional varieties.

\subsection{The Section Conjecture in \'etale homotopy theory}

While hyperbolic curves are \'etale $\K(\pi, 1)$-varieties, many higher-dimensional varieties are not.
For such varieties, it seems unreasonable to expect $\etfdtlgrp(\blank)$ to be a complete invariant.
Instead, one is naturally led to consider their entire \'etale homotopy type.
One contribution of this article is to provide a clear formulation of the Section Conjecture and its pro-$\ell$ variant in terms of \'etale homotopy theory.
We exclusively work with the \emph{profinite \'etale homotopy type} $\ethtpytype(\blank)$, considered as an object of the $\infty$-category $\catpfspc$ of \emph{profinite anima} (see \S\ref{subsec:Sigma-profinite-anima} and \S\ref{subsec:etale-homotopy-type}).

\begin{nul}
    \label{nul:fundamental-fibre-sequence}
    Similarly to the fundamental exact sequence, there is a \emph{fundamental fibre sequence}
    \[
        \ethtpytype(X_{\algclos{k}}) \to \ethtpytype(X) \to \ethtpytype(k)
    \]
    of étale homotopy types for any qcqs scheme $X/k$, see \cite[Corollary 0.5]{fundamental_fibre_sequence}.
\end{nul}

This suggests to replace the set $\Hom_{\Gal_k}^{\out}(\Gal_k, \etfdtlgrp(X))$ of Grothendieck's Section Conjecture with the set $\htpy{\ethtpytype(k), \ethtpytype(X)}_{{\ethtpytype(k)}}$ of homotopy classes of sections of the above fibration $\ethtpytype(X) \to \ethtpytype(k)$.
See \Cref{def:set-of-etale-sections-in-general} for a precise definition.

In terms of \'etale homotopy, the Section Conjecture thus becomes:

\begin{conjecture*}[Section Conjecture]
    The canonical map 
    \[
        X(k) \to \htpy{\ethtpytype(k), \ethtpytype(X)}_{\ethtpytype(k)}, \quad a \mapsto \sqrbr{a_{*}}
    \]
    induced by the functoriality of $\ethtpytype(\blank)$ is a bijection.
\end{conjecture*}

When $X$ is an étale $\K(\pi, 1)$ (e.g. a hyperbolic curve), this recovers the Section Conjecture in its original formulation through fundamental groups, see \Cref{cor:etale-sections-in-terms-of-homotopy-fixed-points}.
Of course, when working over $k = \rR$ one should again replace $X(\rR)$ with $\htpygrp_{0}\!X(\rR)$.
Moreover, by performing a homotopy-theoretic analogue of pushing out along the maximal pro-$\ell$ quotient on the fundamental fibre sequence \ref{nul:fundamental-fibre-sequence}, we are furthermore able to define a geometrically pro-$\ell$ completed variant $\ppfcompl[(\ell)]{\ethtpytype(X)}$ of $\ethtpytype(X)$ together with a canonical comparison map $c_{(\ell)} \from \ethtpytype(X) \to \ppfcompl[(\ell)]{\ethtpytype(X)}$ over $\ethtpytype(k)$, resulting in:

\begin{conjecture*}[pro-$\ell$ Section Conjecture]
    The composite
    \[
        X(k) \xrightarrow{\ethtpytype} \htpy{\ethtpytype(k), \ethtpytype(X)}_{\ethtpytype(k)} \xrightarrow{\blank \circ c_{(\ell)}} \htpy{\ethtpytype(k), \ppfcompl[(\ell)]{\ethtpytype(X)}}_{\ethtpytype(k)}
    \]
    is a bijection.
\end{conjecture*}

We will see in \Cref{cor:sigma-nilpotent-completion-of-etale-homotopy-type-of-curve} that for hyperbolic curves, this again coincides with the usual pro-$\ell$ Section Conjecture.
Note that over number fields and $p$-adic fields, the pro-$\ell$ Section Conjecture for hyperbolic curves is known to be false by work of Hoshi \cite{HoshiNonGeo}; see also \cite[\S 14.7]{stix_book} and \cite{SaidiNonGeo} for additional examples.

\subsection{Overview of results}
\label{subsec:overview}

In order to give a precise statement of our main results, we introduce the following notation.

\begin{definition}[equivariantly triangulable, \ref{def:equivariantly-triangulable-scheme}]
    \label{def:equivariantly-triangulabel-intro}
    A scheme $X/\! \rR$ of finite type is called \emph{equivariantly triangulable} if its complex analytification $X(\cC)$ admits the structure of a finite-dimensional $\Gal_{\rR}$-CW complex.
\end{definition}

\begin{remark}
    \label{rem:equivariantly-triangulable}
    It is known that $X/\! \rR$ of finite type is equivariantly triangulable in the following cases, see \Cref{prop:varieties-that-have-equivariant-CW-structure}.
    \begin{thmlist}
        \item $X$ is smooth over $\rR$.
        \item $X$ is affine.
        \item $X$ is projective over $\rR$.
    \end{thmlist}
    Hofmann \cite[Cor. 4.7]{Hofmann_2009} showed that $X(\cC)$ admits the structure of a finite-dimensional CW complex whenever $X/\rR$ is separated and of finite type. 
    It is plausible that $X$ is even equivariantly triangulable in this generality, but we do not know of a proof.
\end{remark}

Our first main result is, in first approximation, a higher-dimensional and homotopical generalisation of Wickelgren's pro-$2$ real Section Conjecture.
We obtain it as an application of the Sullivan Conjecture\footnote{Formulated by Sullivan in \cite{sullivan}, and proven by Miller \cite{Miller}, Lannes \cite{Lannes} and Carlsson \cite{Carlsson}.} from equivariant homotopy theory.

\begin{theoremA}[{\ref{cor:real-section-conjecture-nilpotent-case-abstract}}]
    \label{thm:A}
    Let $X$ be any equivariantly triangulable scheme over $\rR$.
    The map\footnote{Here, $\fmap_{\ethtpytype(\rR)}(\ethtpytype(\rR), \ppfcompl[(2)]{\ethtpytype(X)})$ is a $2$-profinite refinement of $\map_{\ethtpytype(\rR)}(\ethtpytype(\rR), \ppfcompl[(2)]{\ethtpytype(X)})$, see \Cref{def:fmap}.}

    \[
      \ppfcompl[2]{X(\rR)} \to \fmap_{\ethtpytype(\rR)}(\ethtpytype(\rR), \ppfcompl[(2)]{\ethtpytype(X)})    \]
    adjoint to the composite
    \[
      X(\rR) \xrightarrow{\ethtpytype} \map_{\ethtpytype(\rR)}(\ethtpytype(\rR), \ethtpytype(X)) \xrightarrow{\blank \circ c_{(2)}} \map_{\ethtpytype(\rR)}(\ethtpytype(\rR), \ppfcompl[(2)]{\ethtpytype(X)})
    \]
    is an equivalence of $2$-profinite anima.
    In particular, as follows by passing to connected components, $X/\! \rR$ satisfies the \emph{pro-$2$ Section Conjecture}.
\end{theoremA}

Our second main result is that \emph{geometrically \'etale nilpotent} (see \Cref{def:geometrically-nilpotent-scheme}) varieties over $\rR$ satisfy the full Section Conjecture.
Examples include geometrically \'etale simply connected varieties, as well as abelian varieties.

\begin{theoremB}[\ref{cor:real-section-conjecture-simply-connected-case-abstract}]
    \label{thm:B}
    Let $X$ be any equivariantly triangulable and geometrically étale nilpotent scheme over $\rR$.
    The canonical map 
    \[
        \ppfcompl[2]{X(\rR)} \to \ppfcompl[2]{\map_{\ethtpytype(\rR)}(\ethtpytype(\rR), \ethtpytype(X))}  
    \]
    is an equivalence of $2$-profinite anima.
    In particular, as follows by passing to connected components, $X/\! \rR$ satisfies the \emph{Section Conjecture}.
\end{theoremB}

We derive Theorem B from Theorem A by proving that the canonical comparison map
\[
  \ppfcompl[2]{\map_{\ethtpytype(\rR)}(\ethtpytype(\rR), \ethtpytype(X))} \to \fmap_{\ethtpytype(\rR)}(\ethtpytype(\rR), \ppfcompl[(2)]{\ethtpytype(X)})
\]
is an equivalence whenever $X$ is geometrically \'etale nilpotent.
To achieve this, we develop the basic theory of \emph{pronilpotent} profinite anima. See \S \ref{subsec:nilpotent-sigma-profinite-anima}.

\begin{remark}
    \label{rem:theorems-A-and-B}
    Theorems A and B showcase the following novel phenomena in anabelian geometry:
    \begin{deflist}
        \item In existing results of anabelian geometry, all the {geometric information} is extracted from $\ethtpygrp_{1}$.
              In this sense, Theorem B is \emph{orthogonal} to all other existing results in anabelian geometry.
        \item Theorems A and B provide the first anabelian results for varieties that are \ldots
              \begin{deflist}
                \item not étale $\K(\pi, 1)$'s, as well as
                \item not required to be normal (let alone smooth).
              \end{deflist}
        \item Theorem~A provides a new proof of the classical (pro-$2$) real Section Conjecture for hyperbolic curves, see \Cref{cor:classical-real-section-conjecture}.
    \end{deflist}
\end{remark}

\subsection{Related work}
\label{subsec:related-work}
Within the realm of anabelian geometry, our main Theorems A and B are most closely related to the work of Schmidt and Stix \cite{SchmidtStix}, since we also consistently employ the étale homotopy type to address higher-dimensional anabelian phenomena, along with the existing body of results on the real Section Conjecture.
It should be noted that among the various proofs of the classical real Section Conjecture available in the literature --- most notably Mochizuki's \cite{mzki03}, which contains the first such proof, as well as those by Stix \cite{Stix10}, Vistoli and Bresciani \cite{Vistoli}, Pál \cite{Pal} and Wickelgren's nilpotent version \cite{Wickelgren} --- the methods of the present paper are most similar to those of Pál and Wickelgren:
Pál's proof is certainly influenced by ideas surrounding the Sullivan Conjecture, and Wickelgren actually derives her result by applying the Sullivan Conjecture, albeit in slightly a different fashion than we do.

We should also mention Quick's paper \cite{Quick_2015}:
We follow his strategy to reinterpret the Section Conjecture in terms of homotopy fixed points.
This reinterpretation is crucial for our use of the Sullivan Conjecture.

\subsection{Linear overview}
\label{subsec:linear-overview}
In \S\ref{sec:Sigma-profinite-homotopy-theory}, we present the main homotopy-theoretic ingredients underlying our approach.
We begin by recalling essential background on $\Sigma$-profinite homotopy theory in \S\ref{subsec:Sigma-profinite-anima}.
Next, in \S\ref{subsec:nilpotent-sigma-profinite-anima}, we develop a concise account of nilpotency in $\Sigma$-profinite homotopy theory, culminating in a homotopical analogue of the well-known result that a profinite group is pronilpotent if and only if it splits as the product of its various maximal pro-$\ell$ quotients (see \Cref{thm:characterisation-nilpotent-sigma-profinite-anima}).

In \S\ref{sec:generalised-section-conjecture}, we shift focus to anabelian geometry.
We introduce the \'etale homotopy type in \S\ref{subsec:etale-homotopy-type}, and in \S\ref{subsec:higher-dimensional-section-conjecture}, use it to formulate the generalised (pro-$\ell$) Section Conjecture alluded to above.
Furthermore, we show that for hyperbolic curves, this formulation recovers the classical (pro-$\ell$) Section Conjecture (see \Cref{cor:sigma-nilpotent-completion-of-etale-homotopy-type-of-curve}).

In the final section, \S\ref{sec:proof-of-the-real-section-conjecture}, we establish Theorems A and B.
In \S\ref{subsec:proof-of-theorem-D} we demonstrate a vanishing result for certain nonabelian cohomology groups that, combined with the main result of \S\ref{subsec:nilpotent-sigma-profinite-anima}, allows us to prove the key homotopical result enabling us to deduce Theorem~B from Theorem~A (see \Cref{thm:ppfcompl-commutes-with-htpy-fixed-points-nilpotent-case}).
Finally, in \S\ref{subsec:proofs-of-thm-a-and-b}, we complete the proofs of Theorems A and B.

\subsubsection*{Notation and Conventions}

\begin{enumerate}
    \item We freely make use of the language of $\infty$-categories as developed by Lurie in \cite{HTT}, \cite{HA} and \cite{SAG}.
    \item We will consider ordinary categories (i.e. $1$-categories) as $\infty$-categories via the nerve construction $\nerve(\blank)$, and usually suppress it from the notation.  
    \item We follow Scholze and Clausen's suggestion to replace the term \enquote{space} by \enquote{anima}. 
      We write $\catspc$ for what Lurie calls \enquote{$\infty$-category of spaces} and refer to it as \enquote{$\infty$-category of anima}.
    \item Given an $\infty$-category $\catC$, we write $\Pro(\catC)$ for the $\infty$-category of pro-objects in $\catC$, i.e. the $\infty$-category obtained by freely adjoining cofiltered limits to $\catC$.
      We follow Deligne's notation and denote formal cofiltered limits by $\prolimit$.
      We refer the reader to (the dual of) \cite[\HTTsubsec{5.3.5}]{HTT} and \cite[\SAGsubsec{A.8.1}]{SAG} for more background on pro-objects.
    \item Given an object $t$ and morphisms $x \to t \ot y$ in an $\infty$-category $\catC$, we often write $\map_{t}(y,x)$ instead of $\map_{\overcat{\catC}{t}}(y, x)$ if $\catC$ is inferred from the context.
      Moreover, we write $\htpy{y, x}_{t} \defeqq \htpygrp_{0}\map_{t}(y, x)$ for the set of homotopy classes of maps $y \to x$ over $t$.        
\end{enumerate}

\begin{acknowledgments*}
  This article is based on the author's Ph.D. thesis \cite{phdthesis}.
  We heartily thank Alexander Schmidt and Jakob Stix for their guidance, many insightful discussions, and comments on a draft of this article.
  The author acknowledges support by Deutsche Forschungsgemeinschaft  (DFG) through the Collaborative Research Centre TRR 326 ``Geometry and Arithmetic of Uniformized Structures'', project number 444845124.
\end{acknowledgments*}


  \section{$\Sigma$-profinite homotopy theory}
  \label{sec:Sigma-profinite-homotopy-theory}
  In this section, we collect the general homotopy-theoretic ingredients we need.

\subsection{$\Sigma$-profinite anima and their homotopy and cohomology groups}
\label{subsec:Sigma-profinite-anima}

In \cite[\SAGapp{E}]{SAG}, Lurie develops the basics of \emph{profinite} homotopy theory, i.e. a theory of homotopy types whose homotopy groups naturally carry a profinite structure.
We give a brief overview, and refer the reader to \cite[\SAGapp{E}]{SAG} for full details.

\begin{nul}
  Throughout this section, let $\Sigma$ denote a nonempty set of prime numbers.
  In the following, we will mostly work with
  \begin{thmlist}
    \item $\Sigma = \uppi$ the set of all primes,
    \item $\Sigma = \{p\}$ the set containing a single prime $p$, and
    \item $\Sigma = p'$ the set containing all primes but $p$.
  \end{thmlist}
  We will usually substitute $\Sigma = \{p\}$ with $p$ in formulas, i.e. we write $\catppfspc[p]$ instead of $\catppfspc[\{p\}]$ etc.  
\end{nul}

\begin{definition}[$\Sigma$-finite groups and anima]
    \label{def:Sigma-finite-anima}
    \leavevmode
    \begin{deflist}
      \item A finite group $G$ is said to be \emph{$\Sigma$-finite} if its cardinality is in the multiplicative closure of $\Sigma$.
        We write $\catgrps[\Sigma] \subset \catgrps$ for the full subcategory spanned by the $\Sigma$-finite groups.
      \item An anima $K$ is said to be \emph{$\Sigma$-finite} if the following conditions are satisfied:
            \begin{deflist}
                \item $K$ is truncated, i.e. there exists $N \in \nN$ such that $\htpygrp_{n}(K, k) = 0$ for all $k \in K$ and all $n \geq N$.
                \item The set $\htpygrp_0(K)$ is finite.
                \item For each point $k \in K$ and each integer $n \geq 1$, the group $\htpygrp_n(K,k)$ is $\Sigma$-finite.
            \end{deflist}
      The full subcategory $\catpfinspc[\Sigma] \subset \catspc$ spanned by the $\Sigma$-finite anima is called the \emph{$\infty$-category of $\Sigma$-finite anima}.
    \end{deflist}
\end{definition}
\newpage 
\begin{definition}[$\Sigma$-profinite groups and anima]
    \label{def:Sigma-profinite-spaces}
    \leavevmode
    \begin{deflist}
    \item The category $\catpfgrps[\Sigma]$ is called the \emph{category of\/ $\Sigma$-profinite groups} (or also \emph{pro-$\Sigma$ groups}).
    \item The $\infty$-category $\catppfspc[\Sigma]$ is called the \emph{$\infty$-category of\/ $\Sigma$-profinite anima}.
    In the case $\Sigma = \uppi$ we just say \emph{profinite anima} instead of $\uppi$-profinite anima.
    \end{deflist}
\end{definition}

\begin{remark}
  Write $\cattopgrp$ for the category of topological groups.
  The unique cofiltered-limit-preserving extension
  \[
    \catpfgrps[\uppi] \to \cattopgrp
  \]
  of the functor $\catgrps[\uppi] \subset \cattopgrp, G \mapsto G^{\disc}$ endowing a finite group with the discrete topology is a fully faithful embedding.
  Its image consists of precisely those topological groups that are totally disconnected, compact, and hausdorff.
\end{remark}

\begin{recollection}
    \label{rec:Sigma-profinite-anima-for-different-Sigma}
    Let $\Sigma$ and $\Sigma'$ be two nonempty sets of primes and $\ell$ some prime number.
    \begin{thmlist}
        \item Any object $E \in \catppfspc[\Sigma]$ can be written as a formal cofiltered limit $E \simeq \prolimit_{\alpha} E_{\alpha}$ of $\Sigma$-finite anima $E_{\alpha}$.
        \item Moreover, if $E = \prolimit_{\alpha} E_{\alpha}$ and $B = \prolimit_{\beta} B_{\beta}$ as above, then
          \[
            \map_{\catppfspc[\Sigma]}(E, B) = \limit_{\beta} \colimit_{\alpha} \map(E_{\alpha}, B_{\beta}) \in \catspc.
          \]
        \item Note that $\Sigma' \subset \Sigma$ implies $\catpfinspc[\Sigma'] \subset \catpfinspc[\Sigma]$ and hence also $\catppfspc[\Sigma'] \subset \catppfspc[\Sigma]$.
        \item Since $\Sigma \subset \uppi$, we have an inclusion $\catppfspc[\Sigma] \subset \catpfspc$.
          This inclusion allows us to apply machinery (like homotopy/cohomology groups, \ldots) developed for profinite anima to arbitrary $\Sigma$-profinite anima.
        \item The $1$-point anima $\terminal$ is $\Sigma$-finite for any choice of $\Sigma$.
              A \emph{point} $e$ of $E \in \catppfspc[\Sigma]$ is a map $e \from \terminal \to E$ in $\catppfspc[\Sigma]$, i.e., if $E \simeq \prolimit_{\alpha} E_{\alpha}$, it is a homotopy-coherent choice of points $e_{\alpha}$ of each $E_{\alpha}$.
              We usually write $e \in E$ instead of $e \from \terminal \to E$.
        \item If $E = \prolimit_{\alpha}E_{\alpha}$ is a $\Sigma$-profinite anima, $e \in E$ is a point and $k \geq 0$, then
        \[
            \htpygrp_k(E, e) \defeqq \prolimit_{\alpha} \htpygrp_k(E_{\alpha}, e_{\alpha}),
        \]
        where $e_{\alpha}$ denotes the point on $E_{\alpha}$ induced by $e$, is the \emph{$k$-th homotopy group (resp. set if $k = 0$) of $E$ at $e$}.
        For $k = 0$, it is a profinite set.
        For $k \geq 1$, it is a pro-$\Sigma$ group (abelian for $k \geq 2$).
        \item Similarly, if $E = \prolimit_{\alpha}E_{\alpha}$, $A$ is an abelian group, and $k \geq 0$, then
        \[
           \HH^k(E,A) \defeqq \colimit_{\alpha} \HH^k(E_{\alpha},A)
        \]
        is called the \emph{$k$-th cohomology group of $E$ with values in $A$}.
    \end{thmlist}
\end{recollection}

A crucial property of profinite anima is that the analogue of Whitehead's theorem holds.

\begin{theorem}[profinite Whitehead theorem, see \SAG{}{E.3.1.6}]
  \label{thm:profinite-whitehead-theorem}
  A map $\varphi \from E \to B$ of profinite anima is an equivalence if and only if for all points $e \in E$ and $k \geq 0$, the induced map
        \[
          \varphi_{\! *} \from \htpygrp_k(E, e) \to \htpygrp_k(B, \varphi(e))
        \]
        is an isomorphism of profinite groups (resp., if $k = 0$, profinite sets).
\end{theorem}

When $\Sigma = \{\ell\}$, there is even a purely cohomological characterisation of equivalences.

\begin{theorem}[cohomological characterisation of $\ell$-profinite anima, see {\cite[3.3.15]{dagXIII}}]
  \label{thm:cohomolical-criterion-ell-profinite-equivalences}
        Let $\ell$ be a prime number.
        The following are equivalent for a map $\varphi \from E \to B$ of $\ell$-profinite anima.
        \begin{thmlist}
          \item The map $\varphi \from E \to B$ is an equivalence of $\ell$-profinite anima.
          \item The map $\varphi \from E \to B$ is an \emph{$\fF_{\! \ell}$-equivalence}, i.e. the induced maps
          \[
             \varphi^{\ast} \from \HH^k(B, \fF_{\! \ell}) \to \HH^k(E, \fF_{\! \ell}), \quad k \geq 0
           \]
           on $\fF_{\!\ell}$-cohomology are isomorphisms.
        \end{thmlist}
\end{theorem}

\begin{recollection}[$\Sigma$-profinite completion of anima]
  \label{rec:profinite-completion}
  Let $\Sigma$ be a nonempty set of primes.
  \leavevmode
  \begin{thmlist}
        \item \label{recitem:inclusions-create-finite-limits} The inclusions $\catpfinspc[\Sigma] \subset \catpifinspc \subset \catspc$ create finite limits.
        \item As a consequence, the inclusions $\catppfspc[\Sigma] \subset \catpfspc \subset \Pro(\catspc)$ admit left adjoints
        \[
            \pfcompl{(\blank)} \from \Pro(\catspc) \to \catpfspc \quad \text{and} \quad \ppfcompl[\Sigma]{(\blank)} \from \catpfspc \to \catppfspc[\Sigma],
        \]
        called \emph{profinite completion} and \emph{$\Sigma$-profinite completion}, respectively.
        \item The inclusion $\catspc \subset \Pro(\catspc)$ is left adjoint to the \emph{materialisation} functor
        \[
          \mat \from \Pro(\catspc) \to \catspc, \quad \prolimit_{\alpha}B_{\alpha} \mapsto \limit_{\alpha} B_{\alpha}.
        \]
        \item Since adjoints compose, we see that the restriction of $\ppfcompl[\Sigma]{(\blank)}$ along $\catspc \subset \Pro(\catspc)$ is left adjoint to the restriction of $\mat$ along $\catppfspc[\Sigma] \subset \Pro(\catspc)$.
        \item Pulling back along the unit map $B \to \ppfcompl[\Sigma]{B}$ induces isomorphisms
        \[
           \HH^{k}(\ppfcompl[\Sigma]{B}, A) \to \HH^{k}(B, A), \quad k \geq 0
         \]
         for any profinite anima $B$ and $\Sigma$-finite abelian group $A$.
         In particular, $\htpygrp_0 B = \htpygrp_0 \ppfcompl[\Sigma]{B}$.
         \item \label{recitem:cohomological-finiteness} For $B$ a $\Sigma$-finite anima, $\HH_{n}(B, \zZ)$ is finite for $n \geq 1$.
         In particular, for $\Sigma$-profinite $B$, $\HH^{n}(B, \qQ)$ and $\HH^n(B, \fF_{\!\ell})$ vanish for $n \geq 1$ and almost all primes $\ell$.
         Hence, by \Cref{thm:cohomolical-criterion-ell-profinite-equivalences}, $\ppfcompl[\ell]{B}$ is discrete for almost all primes $\ell$.
  \end{thmlist}
\end{recollection}

\subsubsection{Formal mapping anima}
\label{subsubsec:formal-mapping-anima}

\begin{definition}[formal mapping anima]
  \label{def:fmap}
  Let $E = \prolimit_{\alpha} E_{\alpha}, E' = \prolimit_{\beta}E'_{\beta}$ and $B \in \Pro(\catspc)$.
  \begin{deflist}
  \item We write
  \[
    \fmap(E', E) \defeqq \prolimit_{\alpha} (\colimit_{\beta} \map(E'_{\beta}, E_{\alpha})) \in \Pro(\catspc)
  \]
  for the \emph{formal mapping (pro-)anima of maps $E' \to E$}.
  \item Given a cospan $E' \xrightarrow{s} B \xleftarrow{t} E$, we similarly write $\fmap_B(E', E)$ for the pullback
  \[
    \begin{tikzcd}
      \fmap_B(E', E) \cartesian \arrow[r] \arrow[d] & \fmap(E', E) \arrow[d, "{t \circ \blank}"] \\
      * \arrow[r, "{s}"'] & \fmap(E', B)
    \end{tikzcd}
  \]
  in $\Pro(\catspc)$.
  \end{deflist}
\end{definition}

\begin{nul}
  \label{rem:materialization-of-fmap}
  Note that $\mat[\fmap(E', E)] \simeq \map(E', E)$ and $\mat[\fmap_B(E', E)] \simeq \map_B(E', E)$.
\end{nul}

We will make use of the following slight generalisation of \cite[Proposition 3.3.4]{Anel}.

\begin{lemma}
  \label{lem:cartesian-closed}
  Let $\Sigma$ be a nonempty set of primes.
  If $E$ and $B$ are $\Sigma$-finite anima, so is $\map(E, B)$.
\end{lemma}

Anel states \Cref{lem:cartesian-closed} only in the case $\Sigma = \uppi$, but his proof applies \emph{verbatim} in this more general setting.

\begin{corollary}
  \label{cor:fmap-sometimes-again-Sigma-profinite}
  If $B$ is $\Sigma$-finite and $E = \prolimit_{\alpha} E_{\alpha}$ is $\Sigma$-profinite (over $B$), then $\fmap(B, E)$ (resp. $\fmap_B(B, E)$) is again $\Sigma$-profinite.
\end{corollary}

\subsubsection{Homotopy fixed points and the Sullivan Conjecture}
\label{subsubsec:homotopy-fixed-points}

The $\infty$-category $\catpfspc$ is symmetric monoidal with respect to the product.
As $G$ is a profinite group (i.e. a group object in profinite sets), it is also a group in $\catpfspc$, hence an associative algebra in $(\catpfspc)^{\times}$.
We may therefore consider the $\infty$-category $\RMod_{G}(\catpfspc)$ of right modules over $G$, see \cite[\HAsec{4.2}]{HA} for more details.

\begin{definition}
    \label{def:profinite-spaces-with-G-action}
    Let $G$ be a profinite group.
    We call $\catpfspc[G] \defeqq \RMod_{G}(\catpfspc)$ the \emph{$\infty$-category of profinite anima with (continuous) $G$-action}.
\end{definition}

\begin{remark}
    \label{rem:G-profinite-spaces}
    \leavevmode
    \begin{thmlist}
      \item We usually write $\map_{G}$ instead of $\map_{\catpfspc[G]}$.
      \item If $G$ is a finite group, the above definition unwinds to $\catpfspc[G] = \Fun(\B{G}, \catpfspc)$.
    \end{thmlist}
\end{remark}

The following result is what lets us relate the Section Conjecture with the Sullivan Conjecture.

\begin{theorem}[{\cite[E.6.5.1, E.6.4.4]{SAG}}]
    \label{thm:equivalence-profinite-spaces-with-G-action}
    Let $G$ be a profinite group.
    The construction 
    \[
        G \acts K \quad \mapsto \quad K{\! \modmod \! G} \to \terminal{\! \modmod \! G} = \B{G},
    \] 
    where $K{\! \modmod \! G}$ denotes the \emph{homotopy quotient} of $K$ under the action of $G$, determines an equivalence of $\infty$-categories $\catpfspc[G] \to \overcat{\catpfspc}{\B{G}}$.
    The inverse of this equivalence is given by taking fibres.
\end{theorem}

We now come to the definition of homotopy fixed points.

\begin{definition}
    \label{def:homotopy-fixed-points}
    Let $G$ be a \emph{finite} group and $K \in \catpfspc[G]$ a profinite anima with $G$-action.
    The \emph{(profinite) homotopy fixed point anima} of $K$ is given by $K^{\h{G}} \defeqq \lim_{\B{G}} K \in \catpfspc$.
\end{definition}

A central feature for us is the existence of a homotopy descent spectral sequence converging to the homotopy groups of homotopy fixed point anima:

\begin{proposition}[Quick, {\cite[Theorem 2.16]{Quick_2011}}]
    \label{prop:profinite-homotopy-fixed-point-spectral-sequence}
    Let $G$ be a finite group and $K$ a connected profinite anima with $G$-action.
    Assume $K^{\h{G}} \neq \emptyset$ and let $x \in K^{\h{G}}$.
    Then there is a conditionally convergent homotopy descent spectral sequence
    \[
        \ssE_{2}^{s,t} = \HH^{s}(G, \pfhtpygrp_{t}(K)) \ \converges \ \pfhtpygrp_{t-s}(K^{\h{G}}\! , x), 
    \]
    where the action of $G$ on $\pfhtpygrp_{t}(K)$ depends on the choice of $x$.
\end{proposition}

\begin{proof}
    Denote Quick's model category of profinite spaces by $\catquickpfspc$ (see \cite[]{Quick_2008} for more details).
    We write $\mathcal{M}_{\infty}$ for the underlying $\infty$-category of a model category $\mathcal{M}$, i.e. $\mathcal{M}_{\infty}$ denotes the $\infty$-categorical localisation of the nerve of (the full subcategory of cofibrant objects of) $\mathcal{M}$ with respect to the weak equivalences, see \HA{}{1.3.4.15}.
    By \cite[Corollary 7.4.9]{Harpaz}, there is an equivalence of $\infty$-categories $\catpfspc \simeq \catquickpfspc_{\! \infty}$.
    By \cite[Corollary 2.11]{Quick_2011}, there is a Quillen adjunction between Quick's category of profinite $G$-spaces $\catquickpfspc[G]$ and the slice model category $\overcat{\catquickpfspc\! }{\! \B{G}}$.
    Combining all of these results and using \Cref{thm:equivalence-profinite-spaces-with-G-action}, we obtain a chain of equivalences of $\infty$-categories 
    \[
            \catpfspc[G] \simeq \overcat{(\catpfspc)}{\B{G}}
                         \simeq \overcat{(\catquickpfspc_{\! \infty})\! }{\! \B{G}} 
                         \simeq (\overcat{\catquickpfspc\! }{\! \B{G}})_{\infty} 
                         \simeq \catquickpfspc[G]_{\infty},
    \]
    where, in virtue of \cite[Corollary 7.6.13]{Cisinski}, the second to last equivalence holds since $\B{G}$ is fibrant in $\catquickpfspc$.
    The result hence follows from Quick's profinite homotopy fixed point spectral sequence \cite[Theorem 2.16]{Quick_2011}.
\end{proof}

\begin{remark}
    \label{rem:convergence-of-homotopical-spectral-sequences}
    \leavevmode
    \begin{thmlist}
        \item The differentials in the above spectral sequence are of the form 
        \[
            d_{r} \from \ssE_{r}^{s, t} \to \ssE_{r}^{s+r\!,t+r-1}.
        \]
        Moreover, $\ssE_{r}^{s, t} = \terminal$ if $t - s < 0$ and the diagonal $\ssE_{r}^{s, s}$ is only receiving differentials.
        \item Also note that \cite[Theorem 2.16]{Quick_2011} assumes that $K$ is pointed.
              This assumption is unnecessary: In \cite{BK_ss}, Bousfield explains how to obtain the spectral sequence of an \emph{unpointed} cosimplicial space, as long as one chooses a point of its totalization ($\corresponds K^{\h{G}}$ above).
              So using the machinery of \cite{BK_ss} instead of \cite[\S X.6]{BK}, one obtains the above unpointed variant.
        \item Quick puts unnecessary assumptions on $K$ (resp. $G$) in order to enforce strong convergence of the above spectral sequence (as he later noted himself in \cite[Remark 4.3]{Quick_2015}).
              We refer the reader to \cite[\S 4]{BK_ss} for a detailed discussion of \emph{convergence} of homotopy spectral sequences.
    \end{thmlist}
\end{remark}

To conclude this section, we recall the version of the Sullivan Conjecture we use, due to Lurie, as proven in his course \cite[Lecture 30: Theorem 4]{Lurie_SC}.

\begin{theorem}[Sullivan Conjecture]
    \label{thm:sullivan-conjecture}
    Let $G$ be a finite $p$-group and $K$ a finite-dimensional $G$-CW complex.
    Then the composite map
    \[
        \ppfcompl{(K^{G})} \to \ppfcompl{(K^{\h{G}})} \to (\ppfcompl{K})^{\h{G}}  
    \]
    is an equivalence of $p$-profinite anima.
\end{theorem}

\begin{remark}
    \label{rem:second-map}
    \leavevmode 
    \begin{thmlist}
        \item The map $\ppfcompl{(K^{G})} \to \ppfcompl{(K^{\h{G}})}$ is obtained by applying $\ppfcompl{(\blank)}$ to the canonical map $K^{G} \to K^{\h{G}}$ comparing the ($1$-categorical) limit to the homotopy limit.
        \item The map $\ppfcompl{(K^{\h{G}})} \to (\ppfcompl{K})^{\h{G}}$ is obtained as follows:
              The unit map $K \to \mat[\ppfcompl{K}]$ is $G$-equivariant and therefore induces a map $K^{\h{G}} \to \mat[\ppfcompl{K}]^{\h{G}} \simeq \mat[(\ppfcompl{K})^{\h{G}}]$.
              Since $\ppfcompl{(\blank)}$ is left adjoint to $\mat$, this determines a map $\ppfcompl{(K^{\h{G}})} \to (\ppfcompl{K})^{\h{G}}$.
    \end{thmlist}
\end{remark}

\subsection{Pronilpotent $\Sigma$-profinite anima}
\label{subsec:nilpotent-sigma-profinite-anima}

In this section, we provide a concise treatment of the theory of nilpotency in profinite homotopy theory.

\begin{recollection}[nilpotent group action]
    \label{rec:nilpotent-action}
        Let $\pi$ be a group acting on a (possibly nonabelian) group $G$.
        The action of $\pi$ on $G$ is called \emph{nilpotent} if there exists a finite series of normal subgroups 
        \[
            1 = G_{q} \subset G_{q-1} \subset \ldots \subset G_{0} = G  
        \]
        such that:
        \begin{deflist}
            \item $G_{j-1}/G_{j}$ is abelian and a central subgroup of $G/G_{j}$, and
            \item $G_{j}$ is a $\pi$-subgroup of $G$ and the induced action of $\pi$ on $G_{j-1}/G_{j}$ is trivial.
        \end{deflist}
\end{recollection}

\begin{recollection}[$\pi_{1} \acts \pi_{n}$]
    \label{rec:action-of-fundamental-group-on-higher-homotopy-group}
    Let $K$ be an anima with basepoint $x$.
    Recall that there is a natural action of $\htpygrp_{1}(K, x)$ on the homotopy groups $\htpygrp_{n}(K, x), n \geq 1$.
    Indeed, $\htpygrp_{1}(K, x)$ acts naturally on itself by conjugation.
    One way to see the action on higher homotopy groups is via the universal cover $\trunc_{\geq 1}{\! K}$ of $K$.
    This sits in a fibre sequence 
    \[
        \begin{tikzcd}
            \trunc_{\geq 1}\! K \arrow[r] \arrow[d]
                \arrow[dr, phantom, very near start, "{ \lrcorner }"]
              & K \arrow[d] \\
            \terminal \arrow[r, "{x}"']
              & \B\! \htpygrp_{\leq 1}(K),
        \end{tikzcd}
    \]
    where $\htpygrp_{\leq 1}(K)$ denotes the fundamental groupoid of $K$.
    One has $\htpygrp_{n}(\trunc_{\geq 1}\! K) = \htpygrp_{n}(K, x)$ for $n \geq 2$ and $\htpygrp_{1}(K, x) = \htpygrp_{1}(\B\! \htpygrp_{\leq 1}(K), x)$ naturally acts on the fibre $\trunc_{\geq 1}{\!K}$ (by the classical version of \Cref{thm:equivalence-profinite-spaces-with-G-action}).
    Since $\trunc_{\geq 1}\! K$ is simply connected, this action induces an action on $\htpygrp_{n}(K, x)$ as desired.
\end{recollection}

\begin{definition}[nilpotent anima]
    \label{def:nilpotent-anima}
    Let $K$ be an anima.
    \begin{deflist}
        \item $K$ is said to be \emph{nilpotent} if it is connected and if for any basepoint $x \in K$, the action of $\htpygrp_{1}(K, x)$ on $\htpygrp_{n}(K, x)$ is nilpotent for all $n \geq 1$ (in the sense of \Cref{rec:nilpotent-action}).
        \item $K$ is said to be \emph{componentwise nilpotent} if every connected component $K$ is nilpotent.
        \item We denote by $\catnilspc$ the full subcategory of $\catspc$ spanned by the \emph{componentwise} nilpotent anima.
          Its intersection with $\Ani_{\Sigma}$ is denoted by $\Ani_{\Sigmanil}$.
    \end{deflist}
\end{definition}

\begin{example}
    \label{ex:nilpotent-anima}
    Let $\ell$ be a prime number.
    \begin{thmlist}
        \item \label{rmkitem:sc-is-nilpotent} Any simply connected anima is nilpotent.
        \item More generally, any simple anima (i.e. any connected anima $K$ such that the action of $\htpygrp_{1}(K)$ on $\htpygrp_{n}(K)$ is trivial for all $n \geq 1$) is nilpotent.
        \item \label{rmkitem:p-finite-is-nilpotent} Any $\ell$-finite anima is componentwise nilpotent.
    \end{thmlist}
\end{example}

\begin{definition}[pronilpotent $\Sigma$-profinite anima]
    \label{def:nilpotent-profinite-anima}
    Let $K$ a $\Sigma$-profinite anima.
    \begin{deflist}
        \item $K$ is said to be \emph{pronilpotent} if it can be written as a formal cofiltered limit $K = \prolimit_{\alpha} K_{\alpha}$ with $K_{\alpha}$ componentwise nilpotent and $\Sigma$-finite.
        \item We write $\Pro(\catpfinspc[\Sigma])_{\nil} \simeq \Pro(\catpfinspc[\Sigmanil]) \subset \catppfspc[\Sigma]$ for the full subcategory spanned by the pronilpotent $\Sigma$-profinite anima.
    \end{deflist}
\end{definition}

\begin{lemma}
    \label{lem:simply-connected-profinite-anima-are-nilpotent}
    Let $\ell$ be a prime number.
    \begin{thmlist}
    \item \label{lemitem:sc-implies-pronilpotent} Simply connected profinite anima are pronilpotent.
    \item \label{lemitem:l-profinite-implies-pronilpotent} $\ell$-profinite anima are pronilpotent.
    \end{thmlist}
\end{lemma}

\begin{proof}
    Let $K$ be any simply connected profinite anima.
    By \SAG{}{E.4.6.1}, we can write $K = \prolimit_{\alpha} K_{\alpha}$ with each $K_{\alpha}$ being $\pi$-finite and simply connected.
    But since, by \Cref{ex:nilpotent-anima} \labelcref{rmkitem:sc-is-nilpotent}, simply connected anima are nilpotent, it follows that $K$ is pronilpotent.
    That $\ell$-profinite anima are pronilpotent is an immediate consequence of \Cref{ex:nilpotent-anima} \labelcref{rmkitem:p-finite-is-nilpotent}.
\end{proof}

An important technical advantage of nilpotent anima over simply connected anima are various categorical closure properties shared by the former but not the latter:

\begin{lemma}
    \label{lem:nilpotent-anima-are-stable-under-finite-limits}
    Let $\Sigma$ be a nonempty set of primes.
    \begin{thmlist}
      \item \label{lemitem:inclusion-creates-finite-limits} The inclusion $\catnilspc \subset \catspc$ creates finite limits.
      \item \label{lemitem:inclusion-creates-small-limits} The inclusion $\catppfspc[\Sigma]_{\nil} \subset \catppfspc[\Sigma]$ creates small limits.
            In particular, $\catppfspc[\Sigma]_{\nil}$ has small limits.
    \end{thmlist}
\end{lemma}

\begin{proof}
  \leavevmode
  \begin{thmlist}
    \item It suffices to show that $\catnilspc \subset \catspc$ creates finite products and fibre products.
    That it creates finite products follows immediately from its group-theoretic counterpart.
    Creation of pullbacks is for example proven in \cite[Proposition 4.4.3]{May} for $\class{C} = \catab$ or also \cite[Lemma 7.1]{Dror}.
    \item Since $\catppfspc[\Sigma]_{\nil} = \Pro(\catpfinspc[\Sigmanil])$, it suffices to see that the inclusion
    \[
      \catpfinspc[\Sigmanil] = \catpfinspc[\Sigma] \cap \catpfinspc[\nil] \subset \catpfinspc[\Sigma]
    \]
    creates finite limits.
    This follows immediately from \Cref{rec:profinite-completion} \labelcref{recitem:inclusions-create-finite-limits} and (1). \qedhere
\end{thmlist}
\end{proof}

\subsubsection{A characterisation of pronilpotent \texorpdfstring{$\Sigma$}{Σ}-profinite anima}

Recall that a pro-$\Sigma$ group $G$ is pronilpotent if and only if the canonical homomorphism
\[
  G \to \prod_{\ell \in \Sigma} \pquotient[G][\ell]
\]
of $G$ to the product of its various maximal pro-$\ell$ quotients $\pquotient[G][\ell]$, is an isomorphism.
See e.g. \cite[Prop. 2.3.8]{RZ}.
The goal of this subsection is to prove the following homotopy-theoretic analogue of this fact.

\begin{theorem}[characterisation of pronilpotent $\Sigma$-profinite anima]
    \label{thm:characterisation-nilpotent-sigma-profinite-anima}
    The following are equivalent for a connected $\Sigma$-profinite anima $K$:
    \begin{thmlist}
        \item \label{thmitem:K-is-pronilpotent} $K$ is pronilpotent.
        \item \label{thmitem:map-is-equivalence} The canonical map $K \to \prod_{\ell \! \in \Sigma} \ppfcompl[\ell]{K}$ is an equivalence of $\Sigma$-profinite anima.
    \end{thmlist}
\end{theorem}

We will deduce this by comparing Lurie's $p$-profinite completion $\ppfcompl{(\blank)}$ to the Bousfield-Kan $p$-completion $\BKcompl[p](\blank)$, and using the classical arithmetic fracture square.

\begin{notation}
    \label{not:bousfield-kan-p-completion}
    We write $\BKcompl[p]{(\blank)}$ for the \emph{Bousfield-Kan mod $p$} completion functor, i.e. what Bousfield and Kan in \cite{BK} denote by $R_{\infty}(\blank)$ for $R = \fF_{\! p}$.\!\footnote{Note that, in \cite{BK}, they use the notation $\zZ_{\! p}$ to denote $\fF_{\! p}$.}
\end{notation}

\begin{lemma}
    \label{lem:p-completion-of-pi-finite-nilpotent-anima}
    Let $K$ be a connected $\uppi$-finite anima and $p$ a prime number.
    \begin{thmlist}
        \item The Bousfield-Kan $p$-completion $\BKcompl[p]{K}$ is $p$-finite.
        \item The canonical map $\ppfcompl[p]{K} \to \BKcompl[p]{K}$ induced by $K \to \BKcompl[p]{K}$ is an equivalence.
    \end{thmlist}
    In particular, the $p$-profinite completion $\ppfcompl[p]{K}$ of $K$ is $p$-finite.
\end{lemma}

\begin{proof}
    \leavevmode
    \begin{thmlist}
        \item This is shown in \cite[VII.4.3 (i)]{BK}.
        \item By $(1)$ and \Cref{thm:cohomolical-criterion-ell-profinite-equivalences} it is enough to show that $K \to \BKcompl[p]{K}$ induces an isomorphism on $\fF_{\! p}$-cohomology.
          This is an immediate consequence of the universal property of Bousfield-Kan $p$-completion for \enquote{$p$-good} anima \cite[VII.2.1]{BK} and the fact that $\uppi$-finite anima are $p$-good \cite[VII.4.3 (iii)]{BK}. \qedhere
    \end{thmlist}
\end{proof}

The arithmetic fracture square implies the following characterisation of nilpotent $\Sigma$-finite anima.

\begin{lemma}
    \label{lem:characterisation-pi-finite-nilpotent-anima}
    The following are equivalent for a connected $\Sigma$-finite anima $K$:
    \begin{thmlist}
        \item \label{lemitem:K-is-nilpotent} $K$ is nilpotent.
        \item \label{lemitem:map-is-equivalence} The canonical map $K \to \prod_{\ell \! \in \Sigma} \ppfcompl[\ell]{K}$ is an equivalence.
    \end{thmlist}
\end{lemma}

\begin{proof}
    Since $K$ is $\Sigma$-finite, so is $\prod_{\ell \! \in \Sigma} \BKcompl[\ell]{K} = \prod_{\ell \! \in \Sigma} \ppfcompl[\ell]{K}$ as is seen by combining \Cref{rec:profinite-completion} \labelcref{recitem:cohomological-finiteness} and \Cref{lem:p-completion-of-pi-finite-nilpotent-anima}.
    So both $K$ and $\prod_{\ell \! \in \Sigma} \ppfcompl[\ell]{K}$ have vanishing $\qQ$-cohomology.
    Therefore $\labelcref{lemitem:K-is-nilpotent}$ implies $\labelcref{lemitem:map-is-equivalence}$ by the classical arithmetic fracture square \cite[VI.8.1]{BK}.
    Combining \Cref{rec:profinite-completion} \labelcref{recitem:cohomological-finiteness}, \Cref{ex:nilpotent-anima} \labelcref{rmkitem:p-finite-is-nilpotent}, \Cref{lem:nilpotent-anima-are-stable-under-finite-limits} and \Cref{lem:p-completion-of-pi-finite-nilpotent-anima}, we see that also $\labelcref{lemitem:map-is-equivalence}$ implies $\labelcref{lemitem:K-is-nilpotent}$.
\end{proof}

\begin{proof}[Proof of \Cref{thm:characterisation-nilpotent-sigma-profinite-anima}.]
    Write $K = \prolimit_{\alpha} K_{\alpha}$ for connected $\Sigma$-finite anima $K_{\alpha}$.
    \Cref{lem:p-completion-of-pi-finite-nilpotent-anima} implies that $\ppfcompl[\ell]{K} = \prolimit_{\alpha} \ppfcompl[\ell]{(K_{\alpha})}$.

    We first prove that $\labelcref{thmitem:K-is-pronilpotent}$ implies $\labelcref{thmitem:map-is-equivalence}$.
    If $K$ is pronilpotent, we may assume all the $K_{\alpha}$ to be nilpotent.
    Then $K_{\alpha} \simeq \prod_{\ell \! \in \Sigma} \ppfcompl[\ell]{(K_{\alpha})}$ by \Cref{lem:characterisation-pi-finite-nilpotent-anima}, hence
    \[
            K \simeq \prolimit_{\alpha} \prod_{\ell \! \in \Sigma} \ppfcompl[\ell]{(K_{\alpha})} \simeq \prod_{\ell \! \in \Sigma} \prolimit_{\alpha} \ppfcompl[\ell]{(K_{\alpha})} \simeq \prod_{\ell \! \in \Sigma} \ppfcompl[\ell]{K}
    \]
    as claimed.
    We now prove that also $\labelcref{thmitem:map-is-equivalence}$ implies $\labelcref{thmitem:K-is-pronilpotent}$.
    To this end, note that each $\ppfcompl[\ell]{K}$ is pronilpotent by \Cref{lem:simply-connected-profinite-anima-are-nilpotent} \labelcref{lemitem:l-profinite-implies-pronilpotent}.
    Since, by \Cref{lem:nilpotent-anima-are-stable-under-finite-limits} \labelcref{lemitem:inclusion-creates-small-limits}, pronilpotent $\Sigma$-profinite anima are closed under products, we conclude.
\end{proof}


  \section{The Section Conjecture in \'etale homotopy theory}
  \label{sec:generalised-section-conjecture}
  The goal of this section is to provide a clear formulation of the Section Conjecture and its pro-$\ell$ variant in terms of \'etale homotopy theory.

\subsection{The étale homotopy type}
\label{subsec:etale-homotopy-type}

One major technical advantage of working with $\infty$-categories when developing étale homotopy theory is that the étale homotopy type $\ethtpytype(X)$ is uniquely determined in terms of a universal property.
In order to make this precise, let us briefly recall the $\infty$-categorical version of the étale topos of $X$.

First, we recall the definition of an étale sheaf in the language of $\infty$-categories.

\begin{recollection}[étale sheaves, {\cite[\S A.3.3]{SAG}}]
    Let $X$ be a qcqs scheme and denote by $\etsite{X}$ its (classical, $1$-categorical) small étale site.
    Let $\catC$ be any $\infty$-category.
    \begin{deflist}
        \item Recall that a functor $F \from \etsite{X}[\op] \to \catC$ is an \emph{étale sheaf on $X$ with values in $\catC$} if the following two conditions are satisfied:
            \begin{deflist}
                \item The functor $F$ preserves finite products.
                \item Let $f \from U_{0} \twoheadrightarrow Y$ be an étale surjection and let $U_{\bullet}$ be a \v{C}ech nerve of $f$ (see \cite[\HTTsubsec{6.1.2}]{HTT}), regarded as an augmented simplicial object of $\etsite{X}$.
                    Then the composite map 
                    \[
                        \augmstdsimplex \xrightarrow{U_{\bullet}} \etsite{X}[\op] \xrightarrow{F} \catC
                    \]
                    is a limit diagram, i.e. 
                    \[
                        F(Y) \simeq \limit\left(
                              \begin{tikzcd}[cramped, sep=small]
                            F(U_{0})
                                  \arrow[r, shift left = 0.5em]
                                  \arrow[r, leftarrow]
                                  \arrow[r, shift right = 0.5em]
                                &  F(U_{0} \times_{Y} U_{0})
                                    \arrow[r, shift left = 1em]
                                    \arrow[r, leftarrow, shift left = 0.5em]
                                    \arrow[r]
                                    \arrow[r, leftarrow, shift right = 0.5em]
                                    \arrow[r, shift right = 1em]
                                  &  \ldots
                        \end{tikzcd}\right).
                    \]
            \end{deflist}
        \item The \emph{$\infty$-category of $\catC$-valued étale sheaves on $X$}, denoted by $\Sh_{\et}(X, \catC) = \Sh(\etsite{X}, \catC)$, is the full subcategory of $\Fun(\etsite{X}[\op], \catC)$ spanned by the étale sheaves.
    \end{deflist}
\end{recollection}

Plugging in $\catC = \catspc$ in the above, we arrive at the étale $\infty$-topos of $X$.

\begin{definition}[étale $\infty$-topos]
    \label{def:etale-infinity-topos}
    Let $X$ be a qcqs scheme.
    The $\infty$-category $\ettopos{X} \defeqq \Sh_{\et}(X, \catspc)$ of $\catspc$-valued sheaves for the étale topology on $X$ is called the \emph{étale $\infty$-topos of $X$}.
\end{definition}

We of course have an analogue of \emph{global sections} and \emph{constant sheaves}.

\begin{notation}[global sections]
    \label{not:global-sections}
    The global sections functor 
    \[
        \Gamma_{X, *} \defeqq \Gamma_{\et}(X, \blank) \from \ettopos{X} \to \catspc, \quad F \mapsto F(X)
    \]
    admits a left-exact left adjoint $\Gamma^{*}_{X} \from \catspc \to \ettopos{X}$ that carries an anima $K$ to the \emph{constant sheaf on $X$ (with value $K$)}.
    The adjunction $\Gamma^{*}_{X} \ladj \Gamma_{X, *}$ determines a \emph{geometric morphism} $\ettopos{X} \to \catspc$ of $\infty$-topoi (which is essentially unique since $\catspc$ is the terminal $\infty$-topos).
\end{notation}

\begin{nul}
    We refer the reader to \cite[\HTTch{6}]{HTT} and \cite[\SAGsec{A.3}]{SAG} for more on $\infty$-topoi and sheaves.
\end{nul}

\begin{recollection}[shape of an $\infty$-topos]
	\label{rec:shape}
    Write $\catRTop$ for the $\infty$-category of $\infty$-topoi and (right adjoints in) geometric morphisms.
    The \emph{shape} is a left adjoint functor
    \[
      \htpytype \from \catRTop \to \Pro(\catspc)
    \]
    admitting the following explicit description.
    \begin{deflist}
        \item Given an $\infty$-topos $\topos{X}$, the shape $\htpytype(\topos{X}) \in \Pro(\catspc)$ prorepresents the left exact functor 
                \[
                    \catspc \to \catspc, \quad K \mapsto \Gamma_{*}(\topos{X}, \Gamma^{*}_{\topos{X}} K) \period
                \]
        \item Given a geometric morphism of $\infty$-topoi $f_{*} \from \topos{Y} \to \topos{X}$ with unit $u \from \id_{\topos{Y}} \to f_{*} \circ f^{*}$, the induced map $\htpytype(\topos{Y}) \to \htpytype(\topos{X})$ corresponds to the map 
            \[
                \Gamma_{\topos{X}, *} \ u \ \Gamma_{\topos{X}}^{*} \from \Gamma_{\topos{X}, *} \Gamma_{\topos{X}}^{*} \to \Gamma_{\topos{X}, *} f_{*} f^{*} \Gamma_{\topos{X}}^{*} \simeq \Gamma_{\topos{Y}, *} \Gamma_{\topos{Y}}^{*}
            \]
            in $\Pro(\catspc)^{\op} \subset \Fun(\catspc, \catspc)$.
    \end{deflist}
\end{recollection}

\begin{definition}[étale homotopy type]
    \label{def:etale-shape}
    Let $X$ be a qcqs scheme.
    \begin{deflist}
    \item The proanima $\htpytype(\ettopos{X})$ is called the \emph{étale shape of $X$}.
    \item Its profinite completion $\ethtpytype(X) \defeqq \pfcompl{\htpytype(\ettopos{X})}$ is called the \emph{(profinite) étale homotopy type} of $X$.
    \end{deflist}
\end{definition}

\begin{remark}
  \label{rem:universal-property-of-etale-homotopy-type}
  Let $X$ be a qcqs scheme and $K$ some $\pi$-finite anima.
  Unwinding the definitions, we see that $\ethtpytype(X)$ is defined to be the profinite anima prorepresenting the functor
  \[
    \catpifinspc \to \catspc, \quad K \mapsto \Gamma_{*}(X_{\et}; \Gamma^{*}K).
  \]
  This universal property should be interpreted as a homotopical refinement of comparisons
  \begin{deflist}
    \item $\htpygrp_{1}(\ethtpytype(X), \bar{x}) = \etfdtlgrp(X, \bar{x})$, and
    \item $\HH^{*}(\ethtpytype(X), A) = \HH^{*}_{\et}(X, A)$ for (local systems of) finite abelian groups $A$.
  \end{deflist}
  For example, by the above universal property, we have
  \[
      \HH^{n}(\ethtpytype(X), A) = \htpygrp_{0} \map(\ethtpytype(X), \K(A, n)) 
                                 = \htpygrp_{0} \Gamma_{*}(X_{\et}; \Gamma^{*} \K(A, n)) 
                                 = \HH^{n}_{\et}(X, A),
  \]
  and the comparison of fundamental groups is obtained by similarly considering $\B G$ for varying finite groups $G$.
\end{remark}

\begin{remark}[comparison with Artin--Mazur]
    \label{rem:geometrically-unibranch}
    \leavevmode
    \begin{thmlist}
        \item Hoyois proved in \cite[Proposition 5.1]{Hoyois} that the \enquote{protrunction} of $\htpytype(\ettopos{X})$ (i.e. the pro-left adjoint to the inclusion $\catpfinspc[< \infty] \subset \catspc$ of \emph{truncated} anima into all anima) recovers Artin and Mazur's classical construction of the étale homotopy type introduced in \cite[]{AM}.
        \item Consequently, $\ethtpytype(X)$ always recovers the profinitely completed étale homotopy type of Artin and Mazur.
    \end{thmlist}
\end{remark}

\subsection{The Section Conjecture}
\label{subsec:higher-dimensional-section-conjecture}

It is unreasonable to expect Grothendieck's original formulation of the Section Conjecture to be correct for schemes that are not étale $\K(\pi, 1)$, as $\Hom_{\Gal_{k}}^{\out}(\Gal_k, \etfdtlgrp(X))$ does not contain any information whatsoever on the higher étale homotopy of $X/k$.

In this section, we formulate natural variants of the Section Conjecture and its pro-$\ell$ analogue remedying this deficiency.
Moreover, we show that for hyperbolic curves both variants coincide with their classical counterparts. 

\begin{nul}
    Throughout this section, let $X$ be some qcqs scheme over a field $k$ with separable closure $\algclos{k}$ and absolute Galois group $\Gal_{k} \defeqq \Gal(\algclos{k}/k)$.
\end{nul}

Motivated by the fundamental fibre sequence \ref{nul:fundamental-fibre-sequence}, we replace the set of group-theoretic sections with the following:

\begin{definition}[étale sections]
    \label{def:set-of-etale-sections-in-general}
    Let $X$ be a scheme over a field $k$.
    \begin{deflist}
      \item We call $\map_{\ethtpytype(k)}(\ethtpytype(k), \ethtpytype(X))$ the \emph{anima of étale sections of $X/k$}.
      \item We write
      \[
        \htpy{\ethtpytype(k), \ethtpytype(X)}_{\ethtpytype(k)} \defeqq \htpygrp_0 \map_{\ethtpytype(k)}(\ethtpytype(k), \ethtpytype(X)) 
      \]
      for its set of connected components and refer to it as the \emph{set of étale sections of $X/k$}.
    \end{deflist}
\end{definition}

By functoriality of the étale homotopy type, any $k$-rational point $a \in X(k)$ induces a section
\[
  a_{*} \from \ethtpytype(k) \to \ethtpytype(X)
\]
of $\ethtpytype(X) \to \ethtpytype(k)$.
By taking the corresponding homotopy class, we obtain a replacement for the classical Kummer map.

\begin{definition}[Kummer map]
    \label{def:Kummer-map}
    The map 
    \[
        \kappa_{X/k} \from X(k) \to \etSec{X}{k}, \quad a \mapsto \sqrbr{a_{*}},
    \]
    is called the \emph{Kummer map} of $X/k$.
\end{definition}

This leads to the following general formulation of the Section Conjecture.

\begin{conjecture}[Section Conjecture]
    \label{conj:higher-dimensional-section-conjecture}
    The Kummer map
    \[
        \kappa_{X/k} \from X(k) \to \etSec{X}{k} 
    \]
    is a bijection of sets.
\end{conjecture}

\subsubsection{The pro-$\ell$ Section Conjecture in \'etale homotopy theory}
\label{subsubsec:generalised-pro-ell-Section-Conjecture}

Let $\ell$ be a prime number.

\begin{construction}
  \label{cons:geometrically-pro-ell-homotopy-type}
  Let $X$ be a qcqs scheme over a field $k$.
  \begin{thmlist}
    \item  
    Since, by the Künneth theorem, $\ell$-profinite completion is monoidal, the natural $\Gal_k$-action on $\ethtpytype(X_{\algclos{k}})$ induces a $\Gal_{k}$-action on $\ppfcompl[\ell]{\ethtpytype(X_{\algclos{k}})}$.
    We write $\ppfcompl[(\ell)]{\ethtpytype(X)} \defeqq \ppfcompl[\ell]{\ethtpytype(X_{\algclos{k}})}\modmod\! \Gal_k$ for the corresponding homotopy quotient.
    \item
     The canonical map $c_{\ell} \from \ethtpytype(X_{\algclos{k}}) \to \ppfcompl[\ell]{\ethtpytype(X_{\algclos{k}})}$ is $\Gal_k$-equivariant and hence induces a map
     \[
        c_{(\ell)} \defeqq c_{\ell} \modmod \! \Gal_k \from \ethtpytype(X) \to \ppfcompl[(\ell)]{\ethtpytype(X)}.
     \]
    \item 
     Note that, by \ref{nul:fundamental-fibre-sequence} and \Cref{thm:equivalence-profinite-spaces-with-G-action}, we have a map of fibre sequence
     \[
      \begin{tikzcd}
        \ethtpytype(X_{\algclos{k}}) \arrow[r] \arrow[d, "{c_{\ell}}"'] & \ethtpytype(X) \arrow[r] \arrow[d, "{c_{(\ell)}}"] & \ethtpytype(k) \arrow[d, "{=}"] \\
        \ppfcompl[\ell]{\ethtpytype(X_{\algclos{k}})} \arrow[r] & \ppfcompl[(\ell)]{\ethtpytype(X)} \arrow[r] & \ethtpytype(k).
      \end{tikzcd}
     \]
  \end{thmlist}
\end{construction}

\begin{definition}[geometrically pro-$\ell$ homotopy type]
  \label{def:geometrically-pro-ell-homotopy-type}
  Let $X$ be a qcqs scheme over a field $k$.
  \begin{deflist}
  \item We call $\ppfcompl[(\ell)]{\ethtpytype(X)}$ from \Cref{cons:geometrically-pro-ell-homotopy-type} the \emph{geometrically pro-$\ell$ completed} étale homotopy type of $X/k$.
  \item We call $\map_{\ethtpytype(k)}(\ethtpytype(k), \ppfcompl[(\ell)]{\ethtpytype(X)})$ the \emph{anima of geometrically pro-$\ell$ étale sections} of $X/k$.
  \item We write
  \[
    \htpy{\ethtpytype(k), \ppfcompl[(\ell)]{\ethtpytype(X)}}_{\ethtpytype(k)} \defeqq \htpygrp_0 \map_{\ethtpytype(k)}(\ethtpytype(k), \ppfcompl[(\ell)]{\ethtpytype(X)})
  \]
  for its set of connected components and refer to it as the \emph{set of geometrically pro-$\ell$ étale sections of $X/k$}.
  \end{deflist}
\end{definition}

We can now formulate the pro-$\ell$ Section Conjecture.

\begin{conjecture}[pro-$\ell$ Section Conjecture]
  \label{conj:pro-ell-section-conjecture}
  The composition
  \[
    X(k) \xrightarrow{\kappa_{X/k}} \etSec{X}{k} \xrightarrow{c_{(\ell)} \circ \blank} \htpy{\ethtpytype(k), \ppfcompl[(\ell)]{\ethtpytype(X)}}_{\ethtpytype(k)}
  \]
  is a bijection.
\end{conjecture}

\subsubsection{The real Section Conjecture}
\label{subsubsec:real-section-conjecture}

For $k = \rR$ the real numbers, it is natural to again pass to connected components of the real analytification.

\begin{notation}
    \label{not:real-analytification}
    We simply write $X(\rR)$ for the \emph{real analytification} of $X$.
    Similarly, we write $X(\cC)$ for the \emph{complex analytification} of $X$.
\end{notation}

\begin{conjecture}[real Section Conjecture]
    The Kummer map $\kappa_{X/\!\rR}$ induces a bijection
    \[
        \htpygrp_{0} X(\rR) \to \etSec{X}{\rR}.  
    \]
\end{conjecture}

\begin{conjecture}[real pro-$2$ Section Conjecture]
  The composite
  \[
    \htpygrp_0X(\rR) \xrightarrow{\kappa_{X/\!\rR}} \etSec{X}{\rR} \xrightarrow{c_{(2)} \circ \blank} \htpy{\ethtpytype(\rR), \ppfcompl[(2)]{\ethtpytype(X)}}_{\ethtpytype(\rR)}
  \]
  is a bijection.
\end{conjecture}

\subsubsection{Relation to the classical (pro-$\ell$) Section Conjecture}
\label{subsubsec:relation-to-the-classical-section-conjecture}

The goal of this section is to show that the (pro-$\ell$) Section Conjecture in \'etale homotopy theory coincides with Grothendieck's Section Conjecture for \'etale $\K(\pi, 1)$-varieties.
To this end, we have to understand unpointed homotopy classes of sections of classifying anima of profinite groups.

\begin{proposition}
    \label{prop:homotopy-classes-of-profinite-groupoids}
    Let $\begin{tikzcd}[cramped, sep=small]
        \pi' \arrow[r, "p'"] & G & \arrow[l, two heads, "p"'] \pi
    \end{tikzcd}$
    be a cospan of profinite groups with $p$ surjective.
    Then the mapping $\sqrbr{\varphi} \mapsto \sqrbr{\pfhtpygrp_{1}(\varphi)}$ defines a bijection of sets 
    \[
        \htpygrp_{0} \map_{\B{G}}(\B{\pi'}, \B{\pi}) \to \Hom_{G}(\pi', \pi)_{\Delta},  
    \]
    where $\Delta \defeqq \ker(p)$ acts via conjugation.
\end{proposition}

\begin{proof}
    We use the equivalence $\overcat{(\catpfspc)}{\!\B{G}} \simeq (\overcat{\catquickpfspc}{\!\B{G}})_{\infty}$ with Quick's model category that we already saw in \Cref{prop:profinite-homotopy-fixed-point-spectral-sequence}.
    Since $\pi \twoheadrightarrow G$ is surjective, $\B{\pi} \to \B{G}$ is a fibration in $\catquickpfspc$ (see \cite[Corollary 2.25]{Quick_2008}).
    Since every object is cofibrant, we thus have
    \[
            \htpygrp_{0} \map_{\B{G}}(\B{\pi'}, \B{\pi}) = \htpygrp_{0} \Rmap_{\overcat{\catquickpfspc}{\!\B{G}}}(\B{\pi'}, \B{\pi}) = \htpygrp_{0} \map_{\overcat{\catquickpfspc}{\!\B{G}}}(\B{\pi'}, \B{\pi}).
    \]
    Hence two maps $\varphi_0, \varphi_1 \from \! \B{\pi'} \to \B{\pi}$ over $\B{G}$ are homotopic to each other if and only if there exists a homotopy $h \from \! \B{\pi'} \times \Delta^{1} \to \B{\pi}$ over $\B{G}$ such that the diagram 
    \[
        \begin{tikzcd}
            \B{\pi'} \times \crlybr{0} \arrow[rd, "{\varphi_{0}}"] \arrow[d, hook'] & \\
            \B{\pi'} \times \Delta^{1} \arrow[r, "h" description] & \B{\pi} \\
            \B{\pi'} \times \crlybr{1} \arrow[ru, "{\varphi_{1}}"'] \arrow[u, hook] &
        \end{tikzcd}  
    \] 
    commutes.
    Such a homotopy corresponds to a continuous natural transformation $h \from \varphi_{0} \Rightarrow \varphi_{1}$ of maps of profinite $1$-groupoids\footnote{In the sense of \enquote{$1$-groupoid object in the $1$-category of profinite sets}.} such that the whisker fulfills $p \circ h = \id_{p'}$ (observe that the category of profinite $1$-groupoids embeds fully faithfully into $\catquickpfspc$ via the \enquote{usual} nerve construction).
    Write $c_{g} \from G \to G$ for the conjugation by $g$.
    The choice of such a continuous natural transformation boils down to the choice of an element $\gamma \in \pi$ with the following properties.
    \begin{thmlist}
        \item For all $\sigma \in \pi'$ the diagram 
        \[
            \begin{tikzcd}
                * \arrow[d, "{\gamma}"'] \arrow[r, "{\varphi_{0}(\sigma)}"] & * \arrow[d, "{\gamma}"] \\
                * \arrow[r, "{\varphi_{1}(\sigma)}"'] & *
            \end{tikzcd}  
        \]
        commutes, i.e. $\varphi_{0} = c_{\gamma} \circ \ \varphi_{1}$ (naturality of $h$).
        \item $p(\gamma) = 1$, i.e. $\gamma \in \ker(p) = \Delta$ (since $p \circ h = \id_{p'})$.
        \qedhere
    \end{thmlist}
\end{proof}

\begin{corollary}
    \label{cor:etale-sections-in-terms-of-homotopy-fixed-points}
    Let $X$ be a qcqs étale $\K(\pi, 1)$ scheme over $k$.
    Then the mapping $[s] \mapsto [\htpygrp_{1}(s)]$ defines a bijection of sets
    \[
        \etSec{X}{k} \longrightarrow \Hom_{\Gal_k}^{\out}(\Gal_k, \etfdtlgrp(X)).
    \]
    In particular, the Section Conjecture \ref{conj:higher-dimensional-section-conjecture} precisely recovers Grothendieck's Section Conjecture \ref{conj:section-conjecture} in this case.    
\end{corollary}
\newpage
The comparison to the classical pro-$\ell$ Section Conjecture is based on the following result, see e.g. \cite[Prop. 15]{SchmidtRamification}.

\begin{theorem}
    \label{thm:p-completion-of-etale-homotopy-type-of-curves}
    Let $k$ be a field and $C$ a connected, smooth curve over $k$ that is either incomplete or of strictly positive genus.
    Then the following holds.
    \begin{thmlist}
        \item \label{thmitem:K-pi-1} $C$ is an étale $\K(\pi, 1)$, i.e. $\ethtpytype(C) = \B{\ethtpygrp_{1}(C)}$.
        \item \label{thmitem:ell-K-pi-1} If $k$ is separably closed and $\ell$ is any prime number, then also $\ppfcompl[\ell]{\ethtpytype(C)} = \B{\pquotient[\ethtpygrp_{1}(C)][\ell]}$.
    \end{thmlist}
\end{theorem}

\begin{corollary}
    \label{cor:sigma-nilpotent-completion-of-etale-homotopy-type-of-curve}
    Let $k$ be a field, $C$ a connected and smooth curve over $k$ that is either incomplete or of strictly positive genus and $\ell$ a prime.
    \begin{thmlist}
      \item The Section Conjecture \ref{conj:higher-dimensional-section-conjecture} for $C/k$ recovers the usual Section Conjecture.
      \item We have $\ppfcompl[(\ell)]{\ethtpytype(C)} = \B{\pquotient[\ethtpygrp_1(C)][(\ell)]}$, and the pro-$\ell$ Section Conjecture \ref{conj:pro-ell-section-conjecture} for $C/k$ recovers the usual pro-$\ell$ Section Conjecture.
    \end{thmlist}
\end{corollary}

\begin{proof}
  That the Section Conjecture \ref{conj:higher-dimensional-section-conjecture} recovers the usual Section Conjecture follows immediately from combining \Cref{cor:etale-sections-in-terms-of-homotopy-fixed-points} with \Cref{thm:p-completion-of-etale-homotopy-type-of-curves} \labelcref{thmitem:K-pi-1}.
  Using \Cref{thm:p-completion-of-etale-homotopy-type-of-curves} \labelcref{thmitem:ell-K-pi-1}, we see that
  \[
      \ppfcompl[(\ell)]{\ethtpytype(C)} = \ppfcompl[\ell]{\ethtpytype(C)}\modmod \Gal_{\rR} 
                                     = (\B{\pquotient[\etfdtlgrp(C)][\ell]})\modmod \Gal_{\rR} 
                                     = \B{\pquotient[\etfdtlgrp(C)][(\ell)]}
  \]
  as claimed.
  Therefore, the pro-$\ell$ Section Conjecture \ref{conj:pro-ell-section-conjecture} coincides with the classical pro-$\ell$ Section Conjecture in this case by \Cref{prop:homotopy-classes-of-profinite-groupoids}.
\end{proof}


  \section{Proof of the real (pro-$2$) Section Conjecture}
  \label{sec:proof-of-the-real-section-conjecture}
  In this section, we will finally prove the real (pro-$2$) Section Conjecture.

\subsection{A spectral sequence computation}
\label{subsec:proof-of-theorem-D}

The goal of this section is to prove the following result, the key homotopical input we use to derive the full Section Conjecture from the pro-$2$ Section Conjecture.
 
\begin{proposition}
    \label{thm:ppfcompl-commutes-with-htpy-fixed-points-nilpotent-case}
    Let $G$ be a finite $p$-group and $K$ a connected nilpotent profinite anima with $G$-action.
    Assume that $K^{\h{G}} \neq \emptyset$.
    Then the canonical map 
    \[
        \ppfcompl{(K^{\h{G}})} \to (\ppfcompl{K})^{\h{G}}  
    \]
    is an equivalence of $p$-profinite anima.
\end{proposition}

We start with some preparatory results.

\begin{lemma}
    \label{lem:nonabelian-H1-vanishes}
    Let $\Gamma$ be a profinite group acting continuously on another profinite group $N$.
    If the supernatural orders of $\card{\Gamma}$ and $\card{N}$ are coprime, then $\HH^{1}(\Gamma, N) = \ *$.
\end{lemma}

\begin{proof}
  Since $\HH^1(\Gamma, N)$ classifies sections (up to conjugation) of the tautological exact sequence
  \[
    \begin{tikzcd}[cramped, sep = small]
      1 \arrow[r] & N \arrow[r] & N \semidirectprod \Gamma \arrow[r] & \Gamma \arrow[r] & 1, 
    \end{tikzcd}
  \]
  it suffices to see that all complements of $N$ in $N \semidirectprod \Gamma$ are conjugate to each other.
  This is a direct consequence of the profinite Schur--Zassenhaus theorem \cite[Thm 2.3.15]{RZ}.
\end{proof}

\begin{remark}
  All known proofs of the Schur-Zassenhaus theorem crucially use the Feit--Thompson theorem \cite{FeitThompson}, i.e. that every finite group of odd order is solvable.
\end{remark}

\begin{lemma}
    \label{prop:htpy-fixed-points-away-from-p}
    Let $p$ be a prime number and let $G$ be a finite $p$-group acting on a connected $p'$-profinite anima $K$.
    Assume that $K^{\h{G}} \neq \emptyset$ and let $x \in K^{\h{G}}$.
    Then $\pfhtpygrp_{n}(K^{\h{G}}\! , x) = \pfhtpygrp_{n}(K)^{G}$.
    In particular, $K^{\h{G}}$ is connected.
\end{lemma}

\begin{proof}
    The homotopy fixed point spectral sequence \ref{prop:profinite-homotopy-fixed-point-spectral-sequence}
    \[
        \ssE_{2}^{s, t} = \HH^{s}(G, \pfhtpygrp_{t}(K)) \ \converges \ \pfhtpygrp_{t-s}(K^{\h{G}}\! , x)  
    \]
    attached to $G \acts K$ collapses on the $\ssE_{2}$-page.
    Indeed, since by construction $\ssE_{2}^{s, t} = \terminal$ whenever $t - s < 0$, it suffices to discuss the cases $t \geq s \geq 0$.
    Note that $\ssE_{2}^{s, t} = \HH^{s}(G, \pfhtpygrp_{t}(K))$ vanishes whenever $t \geq 2$ and $s \geq 1$, as we are looking at group cohomology of a finite $p$-group with values in an abelian group prime to $p$.
    When $t = 1$, the only potentially non-trivial term is $\ssE_{2}^{1, 1} = \HH^{1}(G, \pfhtpygrp_{1}(K))$, which vanishes by \Cref{lem:nonabelian-H1-vanishes}.
\end{proof}

\begin{proof}[Proof of \Cref{thm:ppfcompl-commutes-with-htpy-fixed-points-nilpotent-case}]
  As $K$ is nilpotent, the canonical map $K \to \prod_{\ell} \ppfcompl[\ell]{K}$ is an equivalence of profinite anima with $G$-action.
  Hence $K^{\h{G}} \simeq (\prod_{\ell} \ppfcompl[\ell]{K})^{\h{G}} \simeq \prod_{\ell} (\ppfcompl[\ell]{K})^{\h{G}}$.
    Given any prime $\ell$, $(\ppfcompl[\ell]{K})^{\h{G}}$ remains $\ell$-profinite as $\catppfspc[\ell] \subset \catpfspc$ is stable under limits.
    Therefore, since $p$-profinite completion preserves finite products (Künneth formula), we have that 
    \[
        \ppfcompl[p]{(K^{\h{G}})} \simeq (\ppfcompl[p]{K})^{\h{G}} \times \ppfcompl[p]{(\prod_{\ell \neq p} (\ppfcompl[\ell]{K})^{\h{G}})}.   
    \]
    By \Cref{prop:htpy-fixed-points-away-from-p}, $\prod_{\ell \neq p} (\ppfcompl[\ell]{K})^{\h{G}}$ is a product of connected profinite anima, hence itself connected.
    Since it is furthermore $p'$-profinite, we see that $\ppfcompl[p]{(\prod_{\ell \neq p} (\ppfcompl[\ell]{K})^{\h{G}})}$ is contractible.
    Thus $\ppfcompl[p]{(K^{\h{G}})} \to (\ppfcompl[p]{K})^{\h{G}}$ is an equivalence as claimed.
\end{proof}

\subsection{Proofs of Theorems A and B}
\label{subsec:proofs-of-thm-a-and-b}
In the proof of Theorem A, we want to apply the Sullivan Conjecture to $X(\cC)$, so we need to make sure that $X(\cC)$ is a finite-dimensional $\Gal_{\rR}$-CW complex.

\begin{definition}[equivariantly triangulable]
    \label{def:equivariantly-triangulable-scheme}
    A scheme $X$ of finite type over $\rR$ is called \emph{equivariantly triangulable} if there exists a finite-dimensional $\Gal_{\rR}$-CW complex $K$ and a $\Gal_{\rR}$-equivariant homeomorphism $K \to X(\cC)$, where $\Gal_{\rR}$ acts on $X(\cC)$ via complex conjugation.
\end{definition}

The following generalisation of the Riemann existence theorem is what lets us relate the Sullivan Conjecture with the Section Conjecture.

\begin{recollection}[generalised Riemann existence, {\cite[Theorem (12.9)]{AM}}]
    \label{rec:riemann-extistence}
    For any scheme $X$ of finite type over $\cC$, there is a canonical equivalence $\pfcompl{X(\cC)} \isomto \ethtpytype(X)$.
\end{recollection}

\begin{remark}
    \label{rmk:generalised-Riemann-existence}
    Artin and Mazur's version of the generalised Riemann existence theorem assumed $X$ to be pointed and connected.
    Using the modern shape-theoretic definition of the étale homotopy type, these assumptions can be dropped (as we did in the formulation of \Cref{rec:riemann-extistence}), see e.g. \cite[Theorem 11.5.3]{Exodromy}, \cite[Theorem 4.12]{Carchedi} and \cite[Theorem 4.3.10]{Chough}.
\end{remark}

The following constitutes the homotopical incarnation of the pro-$2$ real Section Conjecture:

\begin{theorem}[Theorem A]
    \label{thm:real-section-conjecture-abstract}
    Let $X$ be any equivariantly triangulable scheme over $\rR$.
    The map
    \[
      \ppfcompl[2]{X(\rR)} \to \fmap_{\ethtpytype(\rR)}(\ethtpytype(\rR), \ppfcompl[(2)]{\ethtpytype(X)})
    \]
    adjoint to
    \[
      X(\rR) \xrightarrow{\ethtpytype} \map_{\ethtpytype(\rR)}(\ethtpytype(\rR), \ethtpytype(X)) \xrightarrow{c_{(2)} \circ \blank} \map_{\ethtpytype(\rR)}(\ethtpytype(\rR), \ppfcompl[(2)]{\ethtpytype(X)})
    \]
    is an equivalence of $2$-profinite anima.
\end{theorem}

\begin{proof}
  The equivalence of \Cref{thm:equivalence-profinite-spaces-with-G-action} supplies a canonical equivalence
  \[
    \fmap_{\ethtpytype(\rR)}(\ethtpytype(\rR), \ppfcompl[(2)]{\ethtpytype(X)}) \simeq (\ppfcompl[2]{\ethtpytype(X_{\cC})})^{\h\!\Gal_{\rR}}.
  \]
    We will apply the Sullivan Conjecture in the form of \Cref{thm:sullivan-conjecture}.
    First of all, since $X(\rR) = X(\cC)^{\Gal_{\rR}}$, we also have $\ppfcompl[2]{X(\rR)} = \ppfcompl[2]{(X(\cC)^{\Gal_{\rR}})}$.
    Since $X(\cC)$ is a finite-dimensional $\Gal_{\rR}$-CW-complex, \Cref{thm:sullivan-conjecture} is applicable and shows that furthermore
    \[
            \ppfcompl[2]{X(\rR)} \simeq \ppfcompl[2]{(X(\cC)^{\Gal_{\rR}})} \simeq (\ppfcompl[2]{X(\cC)})^{\h{\! \Gal_{\rR}}} \period
    \]
    Using that $\ppfcompl[2]{X(\cC)} \simeq \ppfcompl[2]{(\pfcompl{X(\cC)})} \simeq \ppfcompl[2]{\ethtpytype(X_{\cC})}$ by the generalised Riemann existence theorem (\ref{rec:riemann-extistence}), we thus see
    \[
            \ppfcompl[2]{X(\rR)} \simeq (\ppfcompl[2]{\ethtpytype(X_{\cC})})^{\h{\! \Gal_{\rR}}},
    \]
    which completes the proof.
\end{proof}

We immediately derive the pro-$2$ real Section Conjecture by applying $\htpygrp_{0}$:

\begin{corollary}[pro-$2$ real Section Conjecture]
    \label{cor:real-section-conjecture-nilpotent-case-abstract}
    Let $X$ be any equivariantly triangulable scheme over $\rR$.
    Then the pro-$2$ Section Conjecture holds for $X$, i.e. the canonical map
    \[
        \htpygrp_{0} X(\rR) \to \htpy{\ethtpytype(\rR), \ppfcompl[(2)]{\ethtpytype(X)}}_{\ethtpytype(\rR)}
    \]
    is a bijection.
\end{corollary}

In order to state the most general form of the real Section Conjecture we can prove, we introduce the following notation.

\begin{definition}[geometrically \'etale nilpotent]
    \label{def:geometrically-nilpotent-scheme}
    A qcqs scheme $X/k$ is said to be \emph{geometrically \'etale nilpotent} if its geometric \'etale homotopy type $\ethtpytype(X_{\algclos{k}})$ is pronilpotent (in the sense of \Cref{def:nilpotent-profinite-anima}).
\end{definition}

\begin{remark}
  \label{rem:simply-connected-and-abelian-varieties-are-nilpotent}
    By \Cref{lem:simply-connected-profinite-anima-are-nilpotent}, any geometrically \'etale simply connected scheme is geometrically \'etale nilpotent.
    Moreover, abelian varieties are geometrically \'etale nilpotent, as they are \'etale $\K(\uppi,1)$-varieties with abelian geometric fundamental group.
\end{remark}

\begin{theorem}[Theorem B]
  \label{thm:Theorem-B-abstract-version}
  Let $X$ be any equivariantly triangulable and geometrically \'etale nilpotent scheme over $\rR$.
  Then the canonical map
  \[
    \ppfcompl[2]{X(\rR)} \to \ppfcompl[2]{\map_{\ethtpytype(\rR)}(\ethtpytype(\rR), \ethtpytype(X))}
  \]
  is an equivalence of $2$-profinite anima.
\end{theorem}

\begin{proof}
  By decomposing $X$ into its connected components, we may assume $X$ to be connected.
  In light of \Cref{thm:real-section-conjecture-abstract}, it is enough to show that the map
  \[
    \ppfcompl[2]{\map_{\ethtpytype(\rR)}(\ethtpytype(\rR), \ethtpytype(X))} \to \fmap_{\ethtpytype(\rR)}(\ethtpytype(\rR), \ppfcompl[(2)]{\ethtpytype(X)})
  \]
  adjoint to
  \[
    \map_{\ethtpytype(\rR)}(\ethtpytype(\rR), \ethtpytype(X)) \xrightarrow{c_{(2)}\circ \blank} \map_{\ethtpytype(\rR)}(\ethtpytype(\rR), \ppfcompl[(2)]{\ethtpytype(X)})
  \]
  is an equivalence.
  If $X(\rR)$ is empty, so is $\fmap_{\ethtpytype(\rR)}(\ethtpytype(\rR), \ppfcompl[(2)]{\ethtpytype(X)})$ by \Cref{thm:real-section-conjecture-abstract}, hence the same holds for $\ppfcompl[2]{\map_{\ethtpytype(\rR)}(\ethtpytype(\rR), \ethtpytype(X))}$.
  We may therefore assume $X(\rR) \neq \emptyset$ and thus that $X$ is geometrically connected.
  Using \Cref{thm:equivalence-profinite-spaces-with-G-action}, we equivalently have to show that the comparison map
  \[
    \ppfcompl[2]{(\ethtpytype(X_{\cC})^{\h\!\Gal_{\rR}})} \to (\ppfcompl[2]{\ethtpytype(X_{\cC})})^{\h\!\Gal_{\rR}}
  \]
  is an equivalence.
  This follows from the assumption that $X$ is geometrically étale nilpotent and \Cref{thm:ppfcompl-commutes-with-htpy-fixed-points-nilpotent-case}.
\end{proof}

Since $2$-profinite completion preserves connected components, this implies the real Section Conjecture.

\begin{corollary}[real Section Conjecture]
    \label{cor:real-section-conjecture-simply-connected-case-abstract}
    Let $X$ be any equivariantly triangulable qcqs and geometrically \'etale nilpotent scheme of finite type over $\rR$.
    Then the Section Conjecture holds for $X$, i.e. the canonical map 
    \[
        \htpygrp_{0} X(\rR) \to \htpy{\ethtpytype(\rR), \ethtpytype(X)}_{\ethtpytype(\rR)}
    \]
    is a bijection.
\end{corollary}

\begin{remark}
    \label{rem:dependency-on-Luries-SC}
    Unfortunately, Lurie's version of the Sullivan Conjecture \ref{thm:sullivan-conjecture} has not been published outside of his set of lecture notes.
    Based on the comparison of Sullivans $p$-adic completion and the Bousfield-Kan $p$-completion $\BKcompl[p](\blank)$, we are still able to derive Theorems A and B, using the more classical versions of the Sullivan Conjecture as follows.
    \begin{thmlist}
        \item\label{thmitem:step-1} Carlsson's version of the Sullivan Conjecture, \cite[Theorem B]{Carlsson}, states that the canonical map 
        \[
                \BKcompl[p]{(K^{G})} \to \BKcompl[p]{(K^{\h{G}})} \to (\BKcompl[p]{\! K})^{\h{G}}
        \]  
        is an equivalence of ordinary anima for every finite $p$-group $G$ and finite-dimensional $G$-CW complex $K$.
        Starting the first step of the proof of \Cref{thm:real-section-conjecture-abstract} with this version of the Sullivan Conjecture, yiels an equivalence $\BKcompl[2]{\! X(\rR)} \simeq (\BKcompl[2]{\! X(\cC)})^{\h{\! \Gal_{\rR}}}$.  
        We can thus continue as in \Cref{thm:real-section-conjecture-abstract}, provided we have an equivalence 
        \[
            \BKcompl[2]{X(\cC)} \simeq \mat[\ppfcompl[2]{X(\cC)}],
        \]
        in which case, because $(\blank)^{\h{\! \Gal_{\rR}}}$ commutes with $\mat$ (both are limits), the $\BKcompl[2]$-version of \Cref{thm:real-section-conjecture-abstract} will result in a natural equivalence
        \[
                \BKcompl[2]{(X(\rR))} \simeq \mat[((\ppfcompl[2]{\ethtpytype(X_{\cC})})^{\h{\!\Gal_{\rR}}})] 
        \]
        of ordinary anima.
        This is still strong enough to deduce Theorems A and B since, by \SAG{}{E.3.1.6}, the materialisation functor $\mat$ is conservative.

        \item The functor $\catspc \to \catspc, K \mapsto \mat[\ppfcompl[p]{K}]$ recovers Sullivan's $p$-profinite completion functor.
              Bousfield-Kan $p$-completion and Sullivan's $p$-profinite completion are known to agree for spaces with degreewise finitely generated $\fF_{\! p}$-homology (see \cite[Thm 6.10]{Friedlander}).
              Therefore, since
              \[
                \HH^{*}(X(\cC); \fF_{\! p}) \cong \HH^{*}_{\et}(X_{\cC}; \fF_{\! p}),
              \]
              the desired comparison $\BKcompl[2]{X(\cC)} \simeq \mat[\ppfcompl[2]{X(\cC)}]$ of \labelcref{thmitem:step-1} indeed holds by \cite[\S XIX, Thm 5.1]{SGA4}.
    \end{thmlist}
    We however chose to stick with using Lurie's version of the Sullivan Conjecture \ref{thm:sullivan-conjecture}, since this better reflects the way we actually arrived at the proofs of Theorems A and B.
\end{remark}

\subsubsection{Examples of equivariantly triangulable schemes}

In order to bring the above theorems to life, we of course need to address the question for which schemes $X/\! \rR$ the complex points $X(\cC)$ actually admit the structure of a finite-dimensional $\Gal_{\rR}$-CW complex, i.e. which $X/\! \rR$ are equivariantly triangulable.

The most general results in this direction that we could find in the literature are summarised in the following proposition:

\begin{proposition}
    \label{prop:varieties-that-have-equivariant-CW-structure}
    Let $X$ be a scheme over $\rR$.
    Then $X$ is equivariantly triangulable, provided one of the following conditions hold:
    \begin{thmlist}
        \item $X$ is smooth over $\rR$.
        \item $X(\cC)$ can be obtained as the zero-locus of a finite family of polynomials in $\rR\sqrbr{T_{1}, \ldots, T_{n}}$ for some $n \in \nN$.
              This includes the cases:
              \begin{thmlist}
                \item $X$ is affine and of finite type over $\rR$.
                \item $X$ is projective over $\rR$.
              \end{thmlist}
    \end{thmlist}
\end{proposition}

\begin{proof}
    We prove the statement case by case.
    \begin{thmlist}
        \item If $X$ is smooth, $X(\cC)$ admits the structure of a smooth $\Gal_{\rR}$-manifold, hence admits a finite-dimensional equivariant $\Gal_{\rR}$-triangulation by the main result of \cite{Illman}.
              Note that, by the conventions regarding equivariant simplicial complexes made in \cite[\S 1]{Illman}, all simplicial complexes considered by Illman are finite-dimensional. 
        \item For $X(\cC) \subset \rR^{n}$ cut out by polynomial equations, the result follows from \cite[Theorem 1.3]{ParkSuh1}.
            \begin{thmlist}
                \item If $X$ is affine and of finite type over $\rR$, then $X(\cC) \subset \cC^{n} \cong \rR^{2n}$ for suitable $n \in \nN$.
                \item If $X$ is projective, then $X(\cC) \subset \PP^{n}(\cC)$ for suitable $n \in \nN$.
                    The claim now follows from the embedding 
                    \[
                        \begin{tikzcd}
                            \PP^{n}(\cC) \arrow[r, closed] & \rR^{2n^{2}}, \quad \sqrbr{z_{0}:\ldots : z_{n}} \mapsto \left( \frac{z_{i}\cdot \bar{z}_{j}}{\sum_{k}\abs{z_{k}}^{2}} \right)_{1\leq i, j\leq n}
                        \end{tikzcd}
                    \]
                    of complex projective $n$-space into real Euclidean $2n^{2}$-space as a compact, algebraic subset, originally due to Mannoury \cite{Mannoury} (see also \cite[Prop. 3.4.6]{RAG}). \qedhere
            \end{thmlist}
    \end{thmlist}
\end{proof}

\begin{remark}
    \label{rem:hofmann}
    In his PhD thesis \cite{Hofmann_2009}, Hofmann shows that the complex points of any separated scheme $X$ of finite type over $\rR$ admit the structure of a finite-dimensional CW complex.
    It seems very plausible that one can actually refine this into the structure of a $\Gal_{\rR}$-CW complex, but we do not know of a proof of this.
\end{remark}

We therefore obtain the following concrete incarnations of \Cref{thm:real-section-conjecture-abstract} and \Cref{thm:Theorem-B-abstract-version}, as well as their corollaries:
\begin{theorem}
  \label{thm:real-section-conjecture}
    Let $X$ be any qcqs scheme that is either smooth, or affine and of finite type, or projective over $\rR$.
    The map
    \[
      \ppfcompl[2]{X(\rR)} \to \fmap_{\ethtpytype(\rR)}(\ethtpytype(\rR), \ppfcompl[(2)]{\ethtpytype(X)})
    \]
    adjoint to
    \[
      X(\rR) \xrightarrow{\ethtpytype} \map_{\ethtpytype(\rR)}(\ethtpytype(\rR), \ethtpytype(X)) \xrightarrow{c_{(2)} \circ \blank} \map_{\ethtpytype(\rR)}(\ethtpytype(\rR), \ppfcompl[(2)]{\ethtpytype(X)})
    \]
    is an equivalence of $2$-profinite anima.
    In particular, the real pro-$2$ Section Conjecture holds for $X$, i.e.
    \[
      \htpygrp_{0}X(\rR) \to \htpy{\ethtpytype(\rR), \ppfcompl[(2)]{\ethtpytype(X)}}_{\ethtpytype(\rR)}
    \]
    is a bijection.
\end{theorem}

\begin{theorem}
  \label{thm:Theorem-B-concrete-version}
  Let $X$ be any qcqs scheme that is either smooth, or affine and of finite type, or projective over $\rR$.
  Assume that $X$ is furthermore geometrically étale nilpotent\footnote{e.g. geometrically étale simply connected or an abelian variety, see \Cref{rem:simply-connected-and-abelian-varieties-are-nilpotent}.}.
  Then the canonical map
  \[
    \ppfcompl[2]{X(\rR)} \to \ppfcompl[2]{\map_{\ethtpytype(\rR)}(\ethtpytype(\rR), \ethtpytype(X))}
  \]
  is an equivalence of $2$-profinite anima.
  In particular, the real Section Conjecture holds for $X$, i.e.
  \[
    \htpygrp_{0}X(\rR) \to \htpy{\ethtpytype(\rR), \ethtpytype(X)}_{\ethtpytype(\rR)}
  \]
  is a bijection.
\end{theorem}

\begin{proof}[Proofs of \Cref{thm:real-section-conjecture} and \Cref{thm:Theorem-B-concrete-version}]\leavevmode

  Combine \Cref{thm:real-section-conjecture-abstract} and \Cref{thm:Theorem-B-abstract-version} with \Cref{prop:varieties-that-have-equivariant-CW-structure}.
\end{proof}

Finally, we note that Theorem A also provides a new proof of the classical real Section Conjecture for hyperbolic curves over $\rR$:

\begin{corollary}
    \label{cor:classical-real-section-conjecture}
    Let $X/\! \rR$ be a hyperbolic curve.
    Then $X$ satisfies the Section Conjecture, i.e. the Kummer map 
    \[
        \kappa_{X/\! \rR} \from \htpygrp_{0} X(\rR) \to \Hom_{\Gal_{\rR}}^{\out}(\Gal_{\rR}, \etfdtlgrp(X))
    \]
    is a bijection.
\end{corollary}

\begin{proof}
    By \cite[Lemma 228]{stix_book} it suffices to show that $\Hom_{\Gal_{\rR}}^{\out}(\Gal_{\rR}, \etfdtlgrp(X)) \neq \emptyset$ implies that $X(\rR) \neq \emptyset$.
    But since, by \Cref{cor:etale-sections-in-terms-of-homotopy-fixed-points}, $\Hom_{\Gal_{\rR}}^{\out}(\Gal_{\rR}, \etfdtlgrp(X)) = \etSec{X}{\rR}$, any such $s$ induces,  via postcomposition with $c_{(2)} \from \ethtpytype(X) \to \ppfcompl[(2)]{\ethtpytype(X)}$, some $\ppfcompl[(2)]{s} \in \htpy{\ethtpytype(\rR), \ppfcompl[(2)]{\ethtpytype(X)}}_{\ethtpytype(\rR)}$.
    Thus $X(\rR) \neq \emptyset$ by \Cref{thm:real-section-conjecture}.
\end{proof}

\begin{remark}
  \label{rem:generalised-pro-2-implies-full-SC}
    We currently do not know whether the analogue of \cite[Lemma 228]{stix_book} holds for the (real) Section Conjecture in \'etale homotopy theory.
    In particular, we do not know whether the pro-$2$ (real) Section Conjecture already implies the (real) Section Conjecture in general.
\end{remark}


  \DeclareFieldFormat{labelnumberwidth}{#1}
  \printbibliography

\end{document}